\documentclass{amsart}

\usepackage{amsmath,color,amssymb,graphics,latexsym}

\newcommand{\bmat}{\left[ \begin{array}}
\newcommand{\emat}{\end{array} \right]}
\newcommand{\ignore}[1]{}

\newtheorem*{lemmaB}{Lemma}

\theoremstyle{definition}

\newtheorem{innerconjecture}{Conjecture}
\newenvironment{customconjecture}[1]
{\renewcommand\theinnerconjecture{#1}\innerconjecture}
{\endinnerconjecture}

\theoremstyle{remark}

\numberwithin{equation}{section}

\begin{document}

\title{Hypergeometric Functions II \\ ($q$--analogues)}

\author{Ian G. Macdonald}

\address{56 High Street, Steventon, Oxfordshire OX13 6RS, England} 

\date{\today}


\maketitle

\tableofcontents

\section*{Foreword}
This is the typewritten version of a handwritten manuscript which was
completed by Ian G. Macdonald in 1987 or 1988.
It is the sequel to the manuscript ``Hypergeometric functions I'' (arXiv:1309.4568).
The two manuscripts are very informal working papers,
never intended for formal publication.
Nevertheless, copies of the manuscripts have circulated widely, giving rise
to quite a few citations in the subsequent 25 years. Therefore it seems
justified to make the manuscripts available for the whole mathematical
community. The author kindly gave his permission that typewritten versions be posted on arXiv. These notes were typeset verbatim by Tierney Genoar and Plamen Koev, supported by the San Jose State University Planning Council, National Science Foundation Grant DMS-1016086, and the
Woodward Fund for Applied Mathematics
at San Jose State University. The Woodward Fund is a gift from the estate of Mrs.\ Marie Woodward in memory of her son, Henry Tynham Woodward. He was an alumnus of the Mathematics Department at San Jose State University and worked with research groups at NASA Ames.
\newpage

\section{}
The classical definition (one variable) is
$$_r\Phi_s(a_1,\ldots,a_r;b_1,\ldots,b_s;x;q) = \sum_{n\geq 0}\frac{(a_1;q)_n\ldots (a_r;q)_n}{(b_1;q)_n\ldots(b_s;q)_n}\cdot\frac{x^n}{(q;q)_n}$$
and we have in particular
\begin{equation}
_0\Phi_0 (x;q) = (x;q)^{-1}_{\infty},
\label{eq_1.1}
\end{equation}
\begin{equation}
_1\Phi_0(a;x;q) = (ax;q)_{\infty}/(x;q)_{\infty}.
\label{eq_1.2}
\end{equation}

In several variables the definitions should be such as to give
\begin{equation}
_0\Phi_0(x_1,\ldots,x_n;q,t) = \prod_{i=1}^n(x_i;q)^{-1}_{\infty},
\label{eq_1.3}
\end{equation}
\begin{equation}
_1\Phi_0(a;x_1,\ldots,x_n;q,t) = \prod_{i=1}^n\frac{(ax_i;q)_{\infty}}{(x_i;q)_{\infty}}.
\label{eq_1.4}
\end{equation}

We have (any $u$)
$$\varepsilon_{u,t}(J_{\lambda}(x;q,t)) = \prod_{(i,j)\in\lambda}(t^{i-1}-q^{j-1}u) = t^{n(\lambda)}\prod_{i\geq 1}(t^{1-i}u;q)_{\lambda_i}$$
and we define
\begin{equation}
(u;q,t)_{\lambda} = \prod_{i\geq 1}(ut^{1-i};q)_{\lambda_i}
\label{eq_1.5}
\end{equation}
so that we have
\begin{equation}
\boxed{\varepsilon_{u,t}(J_{\lambda}) = t^{n(\lambda)}(u;q,t)_{\lambda}.}
\label{eq_1.6}
\end{equation}

In particular,
\begin{equation}
\varepsilon_{0,t}(J_{\lambda}) = t^{n(\lambda)}.
\label{eq_1.7}
\end{equation}

More generally, if $\underline{a} = (a_1,\ldots,a_r)$, define
$$(\underline{a};q,t)_{\lambda} = \prod_{i=1}^r(a_i;q,t)_{\lambda}.$$
Let
$$J_{\lambda}^*(x;q,t) = J_{\lambda}(x;q,t)/\langle J_{\lambda}, J_{\lambda}\rangle_{q,t}$$
so that $(J_{\lambda}^*)$ is the basis of $\Lambda_F$ dual to the basis $(J_{\lambda})$.

\begin{flalign}
\mbox{\underline{Definition}.  Let } \underline {a} = (a_1,\ldots,a_r),\; \underline{b} = (b_1,\ldots,b_s). \mbox{ Then we define}&&
\label{eq_1.8}
\end{flalign}
$$_r\Phi_s(\underline{a};\underline{b};x;q,t) = \sum_{\lambda}\frac{(\underline{a};q,t)_{\lambda}}{(\underline{b};q,t)_{\lambda}}t^{n(\lambda)}J_{\lambda}^*(x;q,t),$$
a formal power series with coefficients in $F(\underline{a},\underline{b})$, i.e., an element of $\hat{\Lambda}_{F(a,b)}$.

Here the number of variables $x_i$ may be finite or infinite.

The next definition, however, seems to be relevant only when the number of variables is finite, say $x=(x_1,\ldots,x_n),$ $y=(y_1,\ldots,y_n)$.

Let
\begin{equation}
J_{\lambda}^*(x,y;q,t) = \frac{J_{\lambda}^*(x;q,t)J_{\lambda}^*(y;q,t)}{\varepsilon_{t^n,t}(J_{\lambda}^*)}
\label{eq_1.9}
\end{equation}
the denominator of which is $J_{\lambda}^*(1,t,\ldots,t^{n-1},q,t)$.
$$ = \frac{J_{\lambda}(x)J_{\lambda}^*(y)}{t^{n(\lambda)}(t^n)_{\lambda}}$$

\begin{flalign}
\mbox{\underline{Definition}.   With } \underline{a}, \underline{b} \mbox{ as in \eqref{eq_1.8}, we define}&&
\label{eq_1.10}
\end{flalign}
$$_r\Phi_s(\underline{a};\underline{b};x,y;q,t) = \sum_{\lambda}\frac{(\underline{a};q,t)_{\lambda}}{(\underline{b};q,t)_{\lambda}}t^{n(\lambda)}J_{\lambda}^*(x,y;q,t).$$
Here the sum is over partitions of length $\leq n$.     (hypergeometric kernel)

The relationship between $_r\Phi_s (x,y)$ and $_r\Phi_s(x)$ is given by
\begin{equation}
\varepsilon^{(y)}_{t^n,t}\,_r\Phi_s(x,y) = \,_r\Phi_s(x)
\label{eq_1.11}
\end{equation}

\begin{proof}
This follows from the definitions, since
$$\varepsilon_{t^n,t}^{(y)}J_{\lambda}^*(x,y) = J_{\lambda}^*(x).$$
\end{proof}

Each such $_r\Phi_s(x,y)$ determines a scalar product on $\Lambda_{F,n}$ for which the $P_{\lambda}$ are pairwise orthogonal and
$$\langle P_{\lambda},Q_{\lambda}\rangle = \frac{(\underline{a})_{\lambda}}{(\underline{b})_{\lambda}(t^n)_{\lambda}}$$
for
$$_r\Phi_s = \sum_{\lambda}\frac{(\underline{a})_{\lambda}}{(\underline{b})_{\lambda}}\cdot\frac{P_{\lambda}(x)Q_{\lambda}(y)}{(t^n)_{\lambda}}.\quad\quad \mbox{(better definition)}$$
\newpage

\section{Particular cases}

\begin{equation}
_0\Phi_0(x;q,t) = \prod_i(x_i;q)_{\infty}^{-1}.
\label{eq_2.1}
\end{equation}

\begin{proof}
\begin{eqnarray*}
_0\Phi_0(x;q,t) &=& \sum_{\lambda}t^{n(\lambda)}J_{\lambda}^*(x;q,t) \\
&=& \sum_{\lambda}\varepsilon_{0,t}(J_{\lambda})J_{\lambda}^*\quad \quad\mbox{by \eqref{eq_1.7}} \\
&=& \varepsilon_{0,t}^{(y)}\Pi(x,y;q,t).
\end{eqnarray*}
Since $\varepsilon_{0,t}(P_r) = (1-t^r)^{-1}$  $(r\geq 1)$, the effect of $\varepsilon_{0,t}$ on the $y$-variables is to specialize $y_i\mapsto t^{i-1}$  $(i\geq 1)$, + hence
$$\varepsilon_{0,t}^{(y)}\Pi(x,y;q,t) = \prod_{i,j}\frac{(x_it^j;q)_{\infty}}{(x_it^{j-1};q)_{\infty}} = \prod_i\frac{1}{(x_i,q)_{\infty}}.$$
\end{proof}

\begin{equation}
_1\Phi_0(a;x;q,t) = \prod_i\frac{(ax_i;q)_{\infty}}{(x_i;q)_\infty}
\label{eq_2.2}
\end{equation}

\begin{proof} By \eqref{eq_1.6} and \eqref{eq_1.8},
$$_1\Phi_0(a;x;q,t) = \sum_{\lambda}\varepsilon_{a,t}(J_{\lambda})J_{\lambda}^* = \varepsilon_{a,t}^{(y)}\Pi (x;y;q,t).$$
Now
$$\Pi(x,y;q,t) = \mathrm{exp}\sum_{r\geq1}\frac{1}{r}\cdot\frac{1-t^r}{1-q^r}p_r(x)p_r(y)$$
and $\varepsilon_{a,t}p_r = \frac{1-a^r}{1-t^r}$, so that
\begin{eqnarray*}
\varepsilon_{a,t}^{(y)}\Pi (x,y;q,t) &=& \mathrm{exp}\sum_{r\geq 1}\frac{1}{r}\cdot\frac{1-a^r}{1-q^r}p_r(x) \\
&=& \Pi (x,1;q,a) \\
&=& \prod_i\frac{(ax_i;q)_{\infty}}{(x_i;q)_{\infty}.}
\end{eqnarray*}
\end{proof}

Notice that \eqref{eq_2.1} is the case $a=0$ of \eqref{eq_2.2}, since $(0;q,t)_{\lambda} = 1$ for all partitions $\lambda$. Thus
$$_r\Phi_s(a_1,\ldots,a_{r-1},\ldots,0;b_1,\ldots,b_s,x;q,t) =\, _{r-1}\Phi_s(a_1,\ldots,a_{r-1};b_1,\ldots,b_s;x;q,t)$$
and likewise if one of the $b_i$ is zero.

\begin{equation}
_1\Phi_0(t^n;x,y;q,t) = \Pi(x,y;q,t)
\label{eq_2.3}
\end{equation}

\begin{proof} We have
\begin{eqnarray*}
_1\Phi_0(t^n;x,y;q,t) &=& \sum_{\lambda}t^{n(\lambda)}(t^n)_{\lambda}J_{\lambda}^*(x,y) \\
&=& \sum_{\lambda}J_{\lambda}(x)J_{\lambda}^*(y)
\end{eqnarray*}
(since $t^{n(\lambda)}(t^n)_\lambda = \varepsilon_{t^n,t}(J_{\lambda})$)
whence the result.
\end{proof}

Next consider the scalar product of Ch.VI, \S9:
$$\langle f,g\rangle'_{q,t} = \frac{1}{n!}[f\bar{g}\Delta]_1,$$
where
$$\Delta = \Delta(x;q,t) = \prod_{i\neq j}\frac{(x_ix_j^{-1};q)_{\infty}}{(tx_ix_j^{-1};q)_{\infty}}$$
and $[\;]_1$ denotes the constant term.

We shall normalize this as follows
$$\langle f,g\rangle''_{q,t} = \langle f,g\rangle'_{q,t}/\langle 1,1\rangle'_{q,t}$$
[In fact ($q$--Dyson conjecture)
$$\langle 1,1\rangle'_{q,t} = \frac{(t;q)^n_{\infty}}{(t^n;q)_{\infty}(q;q)_{\infty}^{n-1}}\cdot\prod_{i=1}^{n-1}(1-t^i)^{-1}$$
 --see e.g., Stembridge for a reasonably simple proof.]

\begin{customconjecture}{(C1)}
$\quad\langle P_{\lambda},P_{\lambda}\rangle''_{q,t} = \varepsilon_{t^n,t}(P_{\lambda})/\varepsilon_{qt^{n-1},t}(Q_{\lambda})$
\label{conjC1}
\end{customconjecture}

We associate with this scalar product the power series
$$\Pi''(x,y;q,t) = \sum_{\lambda}u_{\lambda}(x)v_{\lambda}(y),$$
where $(u_{\lambda}), (v_{\lambda})$ are dual bases of $\Lambda_F$ for the scalar product. Taking $u_{\lambda} = P_{\lambda}$ we have
\begin{eqnarray*}
\Pi''(x,y;q,t) &=& \sum_{\lambda}\frac{P_{\lambda}(x)P_{\lambda}(y)}{\langle P_{\lambda},P_{\lambda}\rangle''} \\
&=& \sum_{\lambda}\frac{P_{\lambda}(x)P_{\lambda}(y)\varepsilon_{qt^{n-1},t}(Q_{\lambda})}{\varepsilon_{t^n,t}(P_{\lambda})} \\
&=& \sum_{\lambda}\frac{\varepsilon_{qt^{n-1},t}(J_{\lambda})J_{\lambda}^*(x)J_{\lambda}^*(y)}{\varepsilon_{t^n,t}(J_{\lambda}^*)} \\
(\mbox{because } P_{\lambda}Q_{\lambda} = J_{\lambda}J_{\lambda}^*) \\ 
&=& \sum_{\lambda}t^{n(\lambda)}(qt^{n-1})_{\lambda}J_{\lambda}^*(x,y) \\
&=& _1\Phi_0(qt^{n-1};x,y;q,t).
\end{eqnarray*}

Then we have

\begin{flalign}
\mbox{Conjecture \ref{conjC1}}  \quad \Longleftrightarrow \quad \Pi''(x,y;q,t) =\, _1\Phi_0(qt^{n-1};x,y;q,t),
\label{eq_2.4}
\end{flalign}

From \eqref{eq_2.4} and \eqref{eq_1.11} we obtain (always assuming \ref{conjC1}) that
\begin{eqnarray*}
\varepsilon^{(y)}_{t^n,t}\Pi''(x,y;q,t) &=&\,_1\Phi_0(qt^{n-1};x;q,t) \\
&=& \prod_i\frac{(qt^{n-1}x_i;q)_{\infty}}{(x_i,q)_{\infty}}\quad\mbox{by }\eqref{eq_2.2} \\
&=& \prod_i(x_i;q)^{-1}_{k(n-1)+1}
\end{eqnarray*}
if $t=q^k$:

\begin{flalign}
\mbox{Conjecture } \mbox{\ref{conjC1}}\quad \Longleftrightarrow \quad \varepsilon_{t^n,t}^{(y)}\,\Pi''(x,y;q,t) = \prod_{i=1}^n\frac{(qt^{n-1}x_i;q)_{\infty}}{(x_i;q)_{\infty}}
\label{eq_2.5}
\end{flalign}
\newpage

\section{Selberg integrals}

Set $t=q^k$ provisionally, where $k$ is a positive integer. Define
\begin{equation}
W_{a,b}(x;q,t) = \prod_{i=1}^nx_i^{a-1}(qx_i;q)_{b-1}\cdot\prod_{1\leq i < j\leq n}\prod_{r=0}^{k-1}(x_i-q^rx_j)(x_i-q^{-r}x_j)
\label{eq_3.1}
\end{equation}

Let $C_n$ denote the unit cube $[0,1]^n$ and define
\begin{equation}
I_{a,b}(f) = \int_{C_n}f(x)W_{a,b}(x)\mathrm{d}_qx
\label{eq_3.2}
\end{equation}
\begin{equation}
J_{a,b}(f) = I_{a,b}(f)/I_{a,b}(1).
\label{eq_3.3}
\end{equation}

$\Big[$The multiple $q$-integral is defined as follows: if f is a function on $C_n$, then
\begin{equation}
\int_{C_n}f(x)\mathrm{d}_q(x) = (1-q)^n\sum_{\alpha\in\mathbb{N}^n}q^{|\alpha|}f(q^{\alpha_1},\ldots ,q^{\alpha_n}).
\label{eq_3.4}
\end{equation}

If $f$ vanishes whenever some $x_i$ is equal to $1$, then
$$\int_{C_n}f(x)\mathrm{d}_qx = (1-q)^n\sum_{\alpha\in\mathbb{N}^n}q^{n+|\alpha|}f(q^{\alpha_1+1},\ldots,q^{\alpha_n+1}),$$
i.e.,
\begin{equation}
\int_{C_n}f(x)\mathrm{d}_qx = q^n\int_{C_n}f(qx)\mathrm{d}_qx.\Big]
\label{eq_3.5}
\end{equation}

\begin{customconjecture}{(C2)}
$J_{a,b}(P_{\lambda}) = \varepsilon_{u,t}(P_{\lambda})\varepsilon_{t^n,t}(P_{\lambda})/\varepsilon_{v,t}(P_{\lambda})$
\label{conjC2}
\end{customconjecture}

\noindent where $u=q^at^{n-1}$,  $v=q^{a+b}t^{2n-2}$.

\smallskip
($t=q^k: \quad u=q^{a'}, v=q^{a'+b'}, a'=a+k(n-1), b'= b+k(n-1)$)
\smallskip

Equivalently, by \eqref{eq_1.6}, \ref{conjC2} $\Longleftrightarrow$

\begin{equation}
J_{a,b}(P_{\lambda}) = \frac{(q^at^{n-1})_{\lambda}}{(q^{a+b}t^{2n-2})_{\lambda}}\varepsilon_{t^n,t}(P_\lambda).
\label{conjC2'}
\tag{C2$'$}
\end{equation}

Conjecture \ref{conjC2} implies
\begin{equation}
I_{a,b}(1)^{-1}\int_{C_n}\,_r\Phi_s(\underline{a};\underline{b};x,y)W_{a.b}(x)\mathrm{d}_qx =\,_{r+1}\Phi_{s+1}(\underline{a},u;\underline{b},v;y),
\label{eq_3.6}
\end{equation}
where as above $u=q^at^{n-1}$,  $v=q^{a+b}t^{2n-2}$.

\begin{proof} We have
\begin{eqnarray*}
I_{a,b}(1)^{-1}\int_{C_n}J_{\lambda}^*(x,y)W_{a,b}(x)\mathrm{d}_qx &=& \frac{J_{a,b}(J_{\lambda}^*(x))J_{\lambda}^*(y)}{\varepsilon_{t^n,t}(J_{\lambda}^*)} \\
&=& \frac{(u)_{\lambda}}{(v)_{\lambda}}J_{\lambda}^*(y)\quad\quad \mbox{by \eqref{conjC2'}}.
\end{eqnarray*}
\end{proof}

Thus integration against $W_{a,b}(x;q,t)$ raises both indices by 1.

We can rewrite the Selberg kernel $W_{a,b}$ in terms of $\Delta(x;q,t)$:--

We have
$$\prod_{r=0}^{k-1}(x_i-q^rx_j)(x_i-q^{-r}x_j) = (-1)^kq^{-k(k-1)/2}(x_ix_j)^k(x_ix_j^{-1};q)_k(x_i^{-1}x_j;q)_k$$
and therefore
\begin{equation}
W_{a,b}(x;q,t) = (-1)^{\alpha}q^{-\beta}\prod_{i=1}^nx_j^{a+k(n-1)-1}(qx_i;q)_{b-1}\Delta(x;q,t),
\label{eq_3.7}
\end{equation}
where $\alpha = k\binom{n}{2}$,  $\beta = \binom{k}{2}\binom{n}{2}$.

\smallskip
So if we define
\begin{equation}
\widetilde{W}_{a,b}(x;q,t) = \prod_{i=1}^nx_i^{a-1}(qx_i;q)_{b-1}\Delta(x;q,t)
\label{eq_3.8}
\end{equation}
and

\begin{equation}
\widetilde{I}_{a,b}(f) = \int_{C_n}f(x)\widetilde{W}_{a,b}(x)\mathrm{d}_qx,
\label{eq_3.9}
\end{equation}
\begin{equation}
\widetilde{J}_{a,b}(f) = \widetilde{I}_{a,b}(f)/\widetilde{I}_{a,b}(1),
\label{eq_3.10}
\end{equation}
we have
\begin{equation}
\widetilde{I}_{a,b}(f) = (-1)^{\alpha}q^{-\beta}I_{a-k(n-1),b}(f)
\label{eq_3.11}
\end{equation}
\begin{equation}
\widetilde{J}_{a,b}(f) = J_{a-k(n-1),b}(f)
\label{3.12}
\end{equation}
so that \ref{conjC2} now takes the form
\begin{equation}
\widetilde{J}_{a,b}(P_{\lambda}) = \frac{(q^a)_{\lambda}}{(q^{a+b\,}t^{n-1})_{\lambda}}\varepsilon_{t^n,t}(P_{\lambda}).
\label{conjC2''}
\tag{C2$''$}
\end{equation}
The value of $I_{a,b}(1)$ is in fact
\begin{equation}
I_{a,b}(1) = n!q^{\gamma}\prod_{i=1}^n\frac{\Gamma_q(ik)\Gamma_q(a+(r-i)k)\Gamma_q(b+(r-i)k)}{\Gamma_q(k)\Gamma_q(a+b+(2r-i-1)k)},
\label{eq_3.13}
\end{equation}
where $\gamma = ka\binom{n}{2}+2k^2\binom{n}{3}$.

If we define
\begin{equation}
\Gamma_{q,n}(a') = \prod_{i=1}^n\Gamma_q(a'-k(i-1))
\label{eq_3.14}
\end{equation}
then \eqref{eq_3.13} takes the form
\begin{equation}
I_{a,b}(1) = n!q^{\gamma}\frac{\Gamma_{q,n}(nk)}{\Gamma_q(k)^n}\cdot\frac{\Gamma_{q,n}(a')\Gamma_{q,n}(b')}{\Gamma_{q,n}(a'+b')},
\label{eq_3.13'}
\tag{3.13$'$}
\end{equation}
where $a'=a+k(n-1)$,  $b'=b+k(n-1)$.
\newpage

\section{Gauss summation for $_2\Phi_1$}

\begin{eqnarray}
\label{eq_4.1}
_2\Phi_1(a_1,a_2;b;c,ct^{-1},\ldots,ct^{1-n};q,t) &=& \prod_{i=1}^n\frac{(a_1^{-1}bt^{1-i};q)_\infty(a_2^{-1}bt^{1-i};q)_\infty}{(bt^{1-i};q)_\infty(a_1^{-1}a_2^{-1}bt^{1-i};q)_\infty} \\
&& \left(=\frac{\Gamma_{n,q}(\beta -\alpha_1)\Gamma_{n,q}(\beta -\alpha_2)}{\Gamma_{n,q}(\beta)\Gamma_{n,q}(\beta -\alpha_1-\alpha_2)}\quad\right),
\nonumber
\end{eqnarray}
where $c=b/(a_1a_2)$.

\begin{proof} From \eqref{eq_3.6} we have
$$_2\Phi_1(a_1,a_2;b;y;q,t) = I_{\alpha,\beta}(1)^{-1}\int_{C_n}\,_1\Phi_0(a_2;x,y)W_{\alpha,\beta}(x)\mathrm{d}_qx$$
where $a_1 = q^{\alpha}t^{n-1}$, $b = q^{\alpha+\beta}t^{2n-2}$ and the effect of replacing $y_i$ by $ct^{1-i}$ replaces $_1\Phi_0(a_2;x,y)$ by $$_1\Phi_0(a_2;ct^{1-n}x) = \prod_i(a_2ct^{1-n}x_i;q)_\infty\big/(ct^{1-n}x_i;q)_\infty.$$

Hence the integral becomes
$$\int_{C_n}\prod_{i=1}^nx_i^{\alpha-1}\cdot\prod^n_{i=1}\frac{(qx_i,q)_\infty}{(q^{\beta} x_i,q)_\infty}\cdot\frac{(a_2ct^{1-n}x_i;q)_\infty}{(ct^{1-n}x_i;q)_\infty}\prod_{i< j}\prod_{r=0}^{k-1}(\;)(\;)\mathrm{d}_qx$$
Now $a_2ct^{1-n} = ba_1^{-1}t^{1-n} = q^\beta$, and $ct^{1-n} = a_2^{-1}q^\beta = q^\gamma$ say.

So finally we have
\begin{eqnarray*}
_2\Phi_1(a_1,a_2;b;(ct^{1-i})_{1\leq i\leq n};q,t) &=& \frac{I_{\alpha,\gamma}(1)}{I_{\alpha,\beta}(1)} \\
&=&\frac{\Gamma_{q,n}(\alpha')\Gamma_{q,n}(\gamma')}{\Gamma_{q,n}(\alpha' +\gamma')}\cdot \frac{\Gamma_{q,n}(\alpha' + \beta')}{\Gamma_{q,n}(\alpha')\Gamma_{q,n}(\beta')}
\end{eqnarray*}
by \eqref{eq_3.13'}, where
$\alpha' = \alpha+ k(n-1)$, $\beta' = \beta +k(n-1)$, $\gamma' = \gamma + k(n-1)$
 so that
$$q^{\alpha'}= a_1,\; q^{\beta'}= b/a_1,\; q^{\gamma'}= a_2^{-1}q^{\beta'}= c = b/a_1a_2,\; q^{\alpha' +\gamma'}= b/a_2.$$
Hence we obtain the formula \eqref{eq_4.1}.

\end{proof}

\subsection*{Additional observation}

\begin{eqnarray*}
\Gamma_{n,q}(a) &=& \prod_{i=1}^n\Gamma_q(a-k(i-1)) \\
&=& \prod_{i=1}^n\frac{(q;q)_\infty}{(q^{a-k(i-1)};q)_\infty}\cdot(1-q)^{-a+k(i-1)}
\end{eqnarray*}
$q$--Saalschutz should be
$$_3\Phi_2(\underline{a};\underline{b};(qt^{n-i})_{1\leq i\leq n};q,t) = \ldots,$$
where some $a_i = q^{-N}$ and $a_1a_2a_3t^{n-1}q=b_1b_2.$

\noindent This will $\Rightarrow$ Gauss as $N \to \infty$.
\newpage

\section{Laplace transform}

The $q$--analogue of $e^x$ is
$$\sum\frac{x^n}{(q;q)_n} = (x;q)^{-1}_\infty$$
and so the analogue of $e^{-x}$ is $(x;q)_\infty$. Let $b\to \infty$ in Conjecture \ref{conjC2}, then $v=q^{a+b}t^{2n-2}\to 0$ and so we have
\begin{equation}
\int_{C_n}P_\lambda(x;q,t)\prod_{i=1}^nx_i^{a-1}(qx_i;q)_\infty\prod_{i<j}\prod_{r=0}^{k-1}(x_i-q^rx_j)(x_i-q^{-r}x_j)\mathrm{d}_qx = I_{a,\infty}(1)^{(q^at^{n-1})_\lambda}\cdot\varepsilon_{t^n,t}(P_\lambda).
\label{eq_5.1}
\tag{5.1}
\end{equation}

Since $\Gamma_q(b) = (q;q)_b\big/(1-q)^{b-1}$, it follows that
$$ \lim_{b\to\infty}\frac{\Gamma_q(b+(n-1)k)}{\Gamma_q(a+b+(2n-i-1)k)} = (1-q)^{a+(n-1)k}$$
and hence from \eqref{eq_3.13} that
$$I_{a,\infty}(1) = n!q^\gamma\prod_{i=1}^n\frac{\Gamma_q(ik)\Gamma_q(a+(n-i)k)}{\Gamma_q(k)}\cdot (1-q)^{na+n(n-1)k}.$$
In the product of gammas the exponent of $(1-q)$is
\begin{eqnarray*}
n(k-1) - \sum_{i=1}^n\big((ik-1)+a+(n-i)k-1\big) &=& nk-n-n(nk+a-2) \\
&=& -(na+n(n-1)k) +n
\end{eqnarray*}
so the total exponent is $n$. So we obtain
\begin{equation}
I_{a,\infty}(1) = n!q^{\gamma}(1-q)^n\prod_{i=1}^n\frac{(q,q)_\infty(q^k;q)_\infty}{(q^{ik};q)_\infty(q^{a+(n-i)k};q)_\infty}.
\label{eq_5.2}
\tag{5.2}
\end{equation}

From \eqref{eq_5.1} it follows that
\begin{equation}
I_{\alpha,\infty}(1)^{-1}\int_{C_n}\,_r\Phi_s (\underline{a},\underline{b};x,y)W_{a,\infty}(x)\mathrm{d}_qx =\, _{r+1}\Phi_s(\underline{a},q^\alpha t^{n-1};\underline{b};y),
\label{eq_5.3}
\tag{5.3}
\end{equation}
raising the $r$--index by 1.

In particular, assuming \ref{conjC1} we have by \eqref{eq_2.3} and \eqref{eq_2.4}
\begin{eqnarray*}
I_{1,\infty}(1)^{-1}\int_{C_n}\Pi(x,y;q,t)W_{1,\infty}(x)\mathrm{d}_qx &=& \, _2\Phi_0(t^n,qt^{n-1};y) \\
&=& I_{k,\infty}(1)^{-1}\int_{C_n}\Pi''(x,y;q,t)W_{k,\infty}(x)\mathrm{d}_qx
\end{eqnarray*}
and 
$$\int_{C_n}\Pi(x,y;q,t)W_{1,\infty}(x)\mathrm{d}_qx = \int_{C_n}\prod_{i,j}\frac{(x_iy_j;q)_\infty}{(tx_iy_j;q)_\infty}\prod_i(qx_i;q)_\infty\prod_{i<j}\prod_{r=0}^{k-1}(x_i -q^{\pm r}x_j)\mathrm{d}_qx.$$

Here
\begin{eqnarray*}
\frac{I_{k,\infty}(1)}{I_{1,\infty}(1)} &=&  q^{k(k-1)\binom{n}{2}}\prod_{i=1}^n\frac{(q^{1+(n-i)k};q)_\infty}{(q^{(n-i+1)k};q)_\infty} \\
&=& q^{k(k-1)\binom{n}{2}}\frac{(q;q)_\infty}{(t^n;q)_\infty}\prod_{i=1}^n(1-t^i)^{-1} \\
&=& q^{k(k-1)\binom{n}{2}}\frac{(q;q)_\infty^n}{(t;q)_\infty^n}\langle1,1\rangle'_{q,t}.
\end{eqnarray*}
\newpage

\section*{Hahn polynomials in one variable}

We shall use the notation
$$(x+1)_a = \frac{(x+a)!}{x!}$$
even if $a,x$ are not integers.

Consider the function
\begin{equation}
F_n^{(a,b)}(x;N) = F_n(x) = (x-n+1)_{a+n}\,(N+1-x)_{b+n},
\label{eq_H1}
\tag{H1}
\end{equation}
where $n,N$ are integers such that $0\leq n\leq N$. Then
\begin{equation}
\Delta^nF_n(x) = \sum_{r=0}^n(-1)^r\binom{n}{r}(x-r+1)_{a+n}\,(N+1-x-n+r)_{b+n},
\label{eq_H2}
\tag{H2}
\end{equation}
in which each term is of the form $(x+1)_a(N+1-x)_b$ multiplied by a polynomial in $x$ (and $a$).
\begin{eqnarray*}
\frac{(x-r+1)_{a+n}}{(x+1)_a} &=& \frac{(x-r+a+n)!}{(x+a)!}\cdot \frac{x!}{(x-r)!} \\
&=& (x+a+1)_{n-r}(x-r+1)_{n-r}
\end{eqnarray*}
and likewise
\begin{eqnarray*}
\frac{(N-x+1-n+r)_{b+n}}{(N-x+1)_a} &=& \frac{(N-x+r+b)!}{(N-x+b)!}\cdot\frac{(N-x)!}{(N-x-n+r)!} \\
&=& (N-x+b+1)_r(N-x-n+r+1)_r.
\end{eqnarray*}

Putting $x=0$ in \eqref{eq_H2}, we obtain
\begin{align}
\nonumber
\Delta^nF_n(0) &= (1)_{a+n}(N+1-n)_{b+n} \\
&= (a+n)!(N+b)!/(N-n)!.
\label{eq_H3}
 \tag{H3}
\end{align}

We define the \underline{$n$th Hahn polynomial} with parameters $a,b$ to be
\begin{equation}
G_n^{(a,b)}(x,N) = \frac{(N-n)!}{N!}\cdot \frac{\Delta_x^n\big((x-n+1)_{a+n}(N+1-x)_{b+n}\big)}{(x+1)_a(N+1-x)_b}.
\label{eq_H4}
\tag{H4}
\end{equation}

It has the following properties:--

\begin{flalign}
\tag{1}
\mbox{\underline{Symmetry.}}  && 
\end{flalign}

 \smallskip

 If $g(x) = f(N-x)$ then
\begin{eqnarray*}
\Delta^ng(x) &=& \sum_{r=0}^n(-1)^r\binom{n}{r}g(x+n-r)\\
&=&\sum_{r=0}^n(-1)^r\binom{n}{r}f(N-x-n+r)\\
&=& (-1)^n(\Delta^nf)(N-x-n)
\end{eqnarray*}

Hence
$$G_n^{(a,b)}(N-x;N) = \frac{(N-n)!}{N!}(-1)^n\frac{\Delta^n((N+1-x)_{a+n}(x-n+1)_{b+n})}{(N+1-x)_a(x+1)_b}$$
i.e.,
\begin{equation}
G_n^{(a,b)}(N-x;N)=(-1)^nG_n^{(b,a)}(x;N).
\label{eq_H5}
\tag{H5}
\end{equation}

Explicitly, we have
\begin{equation}
G_n^{(a,b)}(x;N) = \binom{N}{n}^{-1}\sum_{q+r=n}\frac{(-1)^r}{q!r!}(x+a+1)_q(y+b+1)_r(x-r+1)_q(y-q+1)_r,
\label{eq_H6}
\tag{H6}
\end{equation}
where $y=N-x$.

\begin{flalign}
\tag{2}
\mbox{\underline{Leading term.}} &&
\end{flalign}
\begin{flalign}
\tag{H7} 
G_n^{(a,b)}(x;N) \mbox{ is a polynomial in } x \mbox{ of degree } n, \mbox{ with leading coefficient } (-1)^n\binom{N}{n}^{-1}\binom{a+b+2n}{n}.
\end{flalign}

\begin{proof}. When $a,b$ are positive integers, clearly $(x-n+1)_{a+n}(N+1-x)_{b+n}$ is a polynomial in $x$ of degree $a+b+2n$, with leading coefficient $(-1)^{b+n}$. Hence $G_n^{(a,b)}(x;N)$ is a polynomial of degree
$$(a+b+2n) - n - (a+b) = n$$
with leading coefficient
$$(-1)^n\frac{(N-n)!}{N!}(a+b+n+1)_n = (-1)^n\binom{N}{n}^{-1}\binom{a+b+2n}{n}.$$
Since this is true whenever $a,b\in\mathbb{N}$, it holds generally.
\end{proof}

\begin{flalign}
\tag{3}
\underline{\mbox{Values of } G_n^{(a,b)}(x;N)\mbox{ at }x=0\mbox{ and }x=N.}&&
\end{flalign}

From \eqref{eq_H3} we have
$$G_n^{(a,b)}(0;N) = \frac{(N-n)!}{N!}\cdot\frac{(a+n)!}{(N-n)!}\cdot\frac{N!}{a!} = \frac{(a+n)!}{a!},$$
i.e.,
\begin{equation}
G_n^{(a,b)}(0;N) = (a+1)_n
\label{eq_H8}
\tag{H8}
\end{equation}
and hence by symmetry
\begin{equation}
G_n^{(a,b)}(N;N) = (-1)^n(b+1)_n.
\label{eq_H9}
\tag{H9}
\end{equation}

\begin{flalign}
\tag{4}
\underline{\mbox{Orthogonality}}&&
\end{flalign}

\begin{lemmaB}
(``integration by parts'') For any two functions $f,g$, we have
$$\sum_{x=0}^N(\Delta f)(x)g(x)+\sum_{x=0}^N f(x+1)(\Delta g)(x) = [fg]^{N+1}_0.$$
\end{lemmaB}

\end{document}